# Concentration for norms of infinitely divisible vectors with independent components

CHRISTIAN HOUDRÉ[1], PHILIPPE MARCHAL[2] and PATRICIA REYNAUD-BOURET[3]

[1]*School of Mathematics, Georgia Institute of Technology, Atlanta, GA 30332, USA.*
*E-mail: houdre@math.gatech.edu*
[2]*CNRS and DMA, Ecole Normale Supérieure, 45 rue d'Ulm 75230 Paris Cedex 05, France.*
*E-mail: philippe.marchal@ens.fr*
[3]*CNRS and DMA, Ecole Normale Supérieure, 45 rue d'Ulm 75230 Paris Cedex 05, France.*
*E-mail: patricia.reynaud@ens.fr*

We obtain dimension-free concentration inequalities for $\ell^p$-norms, $p \geq 2$, of infinitely divisible random vectors with independent coordinates and finite exponential moments. Besides such norms, the methods and results extend to some other classes of Lipschitz functions.

*Keywords:* concentration; infinitely divisible laws; norms

## 1. Introduction

The goal of the present paper is to investigate the concentration of measure phenomenon for norms of infinitely divisible random vectors with independent coordinates. Let $X \sim ID(\gamma, 0, \nu)$ be an infinitely divisible (ID) vector without Gaussian component in $\mathbb{R}^d$, and with characteristic function

$$\varphi(t) = \mathbb{E}e^{i\langle t, X\rangle} = \exp\left\{i\langle t, \gamma\rangle + \int_{\mathbb{R}^d}(e^{i\langle t, u\rangle} - 1 - i\langle t, u\rangle\mathbf{1}_{\|u\|_2 \leq 1})\nu(du)\right\}, \qquad (1.1)$$

where $t, \gamma \in \mathbb{R}^d$ and where $\nu \not\equiv 0$ (the Lévy measure) is a positive Borel measure on $\mathbb{R}^d$, without atom at the origin and such that $\int_{\mathbb{R}^d}(1 \wedge \|u\|_2^2)\nu(du) < +\infty$ (throughout, $\langle \cdot, \cdot \rangle$ denotes the Euclidean inner product in $\mathbb{R}^d$, while $\|\cdot\|_2$ is the corresponding Euclidean norm). Properties of $X$ can be read from properties of $\nu$. For example, $X$ has independent components if and only if $\nu$ is supported on the axes of $\mathbb{R}^d$, that is,

$$\nu(dx_1, \ldots, dx_d) = \sum_{k=1}^{d} \delta_0(dx_1) \cdots \delta_0(dx_{k-1})\tilde{\nu}_k(dx_k)\delta_0(dx_{k+1}) \cdots \delta_0(dx_d) \qquad (1.2)$$

for some one-dimensional Lévy measures $\tilde{\nu}_k$, the i.i.d. case corresponding to $\tilde{\nu}_k = \tilde{\nu}$, for $k = 1, \ldots, d$. For simplicity of notation, we shall often assume in the sequel that $X$ has i.i.d. components rather than merely independent ones. It is an easy matter, left to the reader, to transform any i.i.d. result obtained below into an "independent" one; often it just involves introducing $\min_{k=1,\ldots,d}$ or $\max_{k=1,\ldots,d}$ in the notation.

In a seminal work, Talagrand [13] proved a concentration inequality for the product of (one-sided) exponential measures (see also Maurey [9]). This inequality was the first to mix two different norms ($\ell^1$ and $\ell^2$), improving upon some aspects of Gaussian concentration. Rewriting it, in functional form, it asserts that if $X$ is a random vector in $\mathbb{R}^d$ with i.i.d. exponential components and if $f$ is a real-valued Lipschitz function on $\mathbb{R}^d$ such that

$$\exists \alpha, \beta > 0, \forall x, y \in \mathbb{R}^d \qquad |f(x) - f(y)| \leq \min(\alpha \|x - y\|_2, \beta \|x - y\|_1),$$

then there exists a universal constant $K > 0$ such that

$$\mathbb{P}(f(X) - m(f(X)) \geq u) \leq \exp\left(-K \min\left(\frac{x^2}{\alpha^2}, \frac{x}{\beta}\right)\right),$$

where $m(f(X))$ is a median of $f(X)$. What is remarkable here is the dimension-free nature of this concentration inequality. For instance, applying it to the Euclidean norm, we note that $\alpha = \beta = 1$ and that the only dependency in the dimension $d$ is through the median itself. This result of Talagrand, which clearly continues to hold for Lipschitz images of the exponential measure, is actually true for any law satisfying a Poincaré inequality (see Bobkov and Ledoux [2]).

We would here like to obtain dimension-free concentration for infinitely divisible vectors with finite exponential moments, of which the exponential measure is a particular case. The need to have finite exponential moments to obtain dimension-free results is clear, in view of Proposition 5.1 of [13] (see also [1]). On the other hand, by a result of Borovkov and Utev [3], one-dimensional laws satisfying a Poincaré inequality must have a nontrivial absolutely continuous component, thus making our results non-vacuous. The class of infinitely divisible laws is quite encompassing and, for example, on $\mathbb{R}^+$, any log-convex density is infinitely divisible. The situation is more subtle as far as one-dimensional log-concave measures (which necessarily satisfy a Poincaré inequality) is concerned and, for instance, the (infinitely divisible) gamma law with parameters $\alpha > 0$ and $t > 0$, and with density $\alpha^t x^{t-1} e^{-\alpha x}/\Gamma(t), x > 0$, is log-concave if and only if $t \geq 1$. Let us also mention that double Wiener–Itô integrals form another important example of infinitely divisible laws (we refer the reader to Sato [12] for a comprehensive introduction to infinitely divisible laws).

For general Lipschitz functions and general ID vectors, generic results have already been obtained, but when specialized to vectors with i.i.d. components, they are not always dimension-free (we refer to [6] for more precise statements). In fact, it is not clear whether or not an extra assumption, such as convexity, might be needed in order to obtain dimension-free concentration for generic Lipschitz functions. As shown below, for $\ell^p$-norms, $p \geq 2$, we do obtain dimension-free concentration.

Let us state a first result for the Euclidean norm.

**Theorem 1.** *Let $X = (X_1, \ldots, X_d)$ be an ID vector with i.i.d. coordinates, characteristic function (1.1) and Lévy measure as in (1.2). Let $\mathbb{E}e^{t\|X\|_2} < +\infty$ for some $t > 0$, let $M = \sup\{t > 0 : \mathbb{E}e^{t|X_1|} < +\infty\}$ and also let $l = -\log \mathbb{E}[e^{-X_1^2}]$. Then, for all $x > 0$,*

$$\mathbb{P}(\|X\|_2 \geq \mathbb{E}\|X\|_2 + x) \leq e^{-\sup_{0 \leq t \leq T}[tx - \int_0^t 2g(s)\,ds]}, \tag{1.3}$$

*where, for $0 < t < M$,*

$$g(t) = \left(8 + \frac{12 \log(2)}{l}\right) \int_{\mathbb{R}} |u|(e^{t|u|} - 1)\tilde{\nu}(du) + \frac{8}{l} \int_{\mathbb{R}} |u|^3 (e^{t|u|} - 1)\tilde{\nu}(du)$$

*and where $T$ is such that for all $t \leq T$, $tg(t) \leq 1/2$.*

Inequality (1.3) does recover Talagrand's inequality for the Euclidean norm (up to the value of the constants). Indeed, for the symmetric exponential law, in which case $\tilde{\nu}(du)/du = e^{-|u|}/|u|, u \in \mathbb{R}, u \neq 0$, we obtain $T \simeq 0.06$. Moreover, since $g(0) = 0$, there exists a $C > 0$ such that $\int_0^t 2g(s)\,ds \leq Ct^2$ for all $t \leq T$. Taking $t = x/(2C)$ for $x \leq 2CT$ and $t = T$ otherwise, we get bounds of the form $\exp(-x^2/(4C))$ for $x \leq 2CT$ and $K \exp(-Tx)$ for $x > 2CT$.

We next obtain a result for general $\ell^p$-norms, $p \geq 2$, but under some assumptions on the law of $X$.

**Theorem 2.** *Let $X$ be as in Theorem 1 and let $X_1$ be either symmetric or non-negative. Let $2 \leq p < \infty$. Then, for all $0 < x < h_p(M^-)$,*

$$\mathbb{P}(\|X\|_p - \mathbb{E}\|X\|_p \geq x) \leq \exp\left(-\int_0^x h_p^{-1}(s)\,ds\right), \tag{1.4}$$

*where the (dimension-free) function $h_p$ is given by*

$$h_p(t) = p^2 \int_{\mathbb{R}} \left[\left(1 + \frac{4^{1/p}|u|}{m_p^{1/p}}\right)^{2p-2} + 2^{2p+1} \frac{m_{2p}}{m_p^2}\right] |u|(e^{t|u|} - 1)\tilde{\nu}(du),$$

*where $0 < t < M$ and where, for any $q > 1$, $m_q = \mathbb{E}[|X_1|^q]$.*

The hypotheses on $X$ may seem restrictive, but we shall see that a similar result (see Theorem 5 in Section 4) holds under far more general conditions. However, the dimension-free bound we obtain in that general framework is more complicated to express.

Note that Theorem 2 also recovers a bound of the form $\exp(-\min(cx^2, c'x))$, as in Talagrand's result. Moreover, the constant $c'$ is now asymptotically optimal. More precisely, suppose that $X$ is infinitely divisible (without Gaussian component), one-dimensional and satisfies

$$\log \mathbb{P}(|X| \geq x) \sim -\lambda_0 x \tag{1.5}$$

as $x \to \infty$, for some constant $\lambda_0 > 0$. Then, for every $\lambda < \lambda_0$,

$$\int_{\mathbb{R}} |u|(\mathrm{e}^{\lambda|u|} - 1)\tilde{\nu}(\mathrm{d}u) < \infty$$

and, therefore, $h_p^{-1}$ is well defined on $[0, \lambda_0)$. It follows that if we study the $\ell^p$-norm of $(X_1, \ldots, X_n)$, where the $X_i$ are i.i.d. and have the same law as $X$, Theorem 2 gives, for every $\varepsilon > 0$, a bound of order $\exp(-(\lambda_0 - \varepsilon)x)$ for large $x$. In view of (1.5), we see that this bound is optimal, up to some subexponential factor.

For instance, suppose that $\tilde{\nu}$ is concentrated on $\mathbb{R}_+$, has a density $k$ and there exist two constants $\lambda_0, q > 0$ such that, as $x \to \infty$,

$$k(x) \asymp x^{-q}\mathrm{e}^{-\lambda_0 x},$$

where, as usual, $\asymp$ indicates that the ratio of the two quantities is bounded, above and below, as $x \to \infty$. Then, if $2 \leq p < q/2$, Theorem 2 gives, for large enough $x$, a bound of the form

$$\mathbb{P}(\|X\|_p - \mathbb{E}\|X\|_p \geq x) \leq c(x)\mathrm{e}^{-\lambda_0 x},$$

where $\log c(x)/\log x \to 0$ as $x \to \infty$. On the other hand, if $p \geq \sup(2, q/2)$, then for every $\lambda < \lambda_0$, if $x$ is large enough,

$$\mathbb{P}(\|X\|_p - \mathbb{E}\|X\|_p \geq x) \leq c(x)\mathrm{e}^{-\lambda x}.$$

Moreover, if the Lévy measure $\nu$ has bounded support, Theorem 2 gives a bound of order $\exp(-x \log x)$ for large $x$. This is known to be the right order of magnitude for a Poisson random variable (see, e.g., [6]). In turn, this kind of bound entails the existence of more-than-exponential moments. The most precise result we obtain is the following dimension-free extension of the results of [6, 11].

**Theorem 3.** *Let $X$ be as in Theorem 1, let $\tilde{\nu}$ have bounded support and let $R = \inf\{\rho : \tilde{\nu}(|x| > \rho) = 0\}$. Then,*

$$\mathbb{E}[\mathrm{e}^{(\|X\|_2/R)\log(\lambda\|X\|_2/R)}] < +\infty,$$

*for all $\lambda$ such that $\lambda V^2/R^2 < 1/\mathrm{e}$, where $V^2 = 8\int_{\mathbb{R}} |u|^2 \tilde{\nu}(\mathrm{d}u)$.*

Further results of a similar flavor, dealing with projections, $\ell^p$-norms or integrals with respect to a Poisson process, are given in the remainder of this paper. All these results are based on a covariance formula that can be derived from a result in [5]. This formula, together with its first applications, is proved in Section 2. Theorem 1 is then proved in Section 3. In Section 4, we state and prove Theorem 5, which is a generalization of Theorem 2. The last section is devoted to the proof of Theorem 3.

## 2. The covariance formula and its first applications

### 2.1. The covariance formula

The result at the root of every proof in this paper is the following one.

**Proposition 1.** *Let $X = (X_1, \ldots, X_d) \sim ID(\gamma, 0, \nu)$ have independent components and be such that $\mathbb{E} e^{t\|X\|_2} < +\infty$ for some $t > 0$. Let $f : \mathbb{R}^d \to \mathbb{R}$ be such that $\mathbb{E} f(X) = 0$ and let there exist $b_k \in \mathbb{R}$, $k = 1, \ldots, d$, such that $|f(x + ue_k) - f(x)| \leq b_k|u|$ for all $u \in \mathbb{R}$, $x \in \mathbb{R}^d$. Let $M = \sup\{t > 0 : \forall k = 1, \ldots, d, \mathbb{E} e^{tb_k|X_k|} < +\infty\}$. Then, for all $0 \leq t < M$,*

$$\mathbb{E} f e^{tf} \leq \int_0^1 \mathbb{E}_z \left[ \sum_{k=1}^d \int_{\mathbb{R}} \frac{|f(U + ue_k) - f(U)|^2 + |f(V + ue_k) - f(V)|^2}{2} \right.$$
$$\left. \times e^{tf(V)} \left( \frac{e^{tb_k|u|} - 1}{b_k|u|} \right) \tilde{\nu}_k(\mathrm{d}u) \right] \mathrm{d}z,$$

*where the expectation $\mathbb{E}_z$ is with respect to the ID vector $(U, V)$ in $\mathbb{R}^{2d}$ of parameter $(\gamma, \gamma)$ and with Lévy measure $z\nu_1 + (1 - z)\nu_0$, $0 \leq z \leq 1$. The measure $\nu_0$ is given by*

$$\nu_0(\mathrm{d}u, \mathrm{d}v) = \nu(\mathrm{d}u)\delta_0(\mathrm{d}v) + \delta_0(\mathrm{d}u)\nu(\mathrm{d}v), \qquad u, v \in \mathbb{R}^d,$$

*while $\nu_1$ is the measure $\nu$ supported on the main diagonal of $\mathbb{R}^{2d}$.*

An important feature of this proposition is the fact that the first marginal of $(U, V)$ is $X$ and so is its second marginal. Therefore, the main problem in estimating the right-hand side of the inequality in Proposition 1 will be to decouple $U$ and $V$, that is, to split the product $|f(U + ue_k) - f(U)|^2 e^{tf(V)}$ without changing the term $e^{tf(V)}$. To do so, a first attempt could be to use a supremum.

**Corollary 1.** *Let $X = (X_1, \ldots, X_d) \sim ID(\gamma, 0, \nu)$ have independent components and be such that $\mathbb{E} e^{t\|X\|_2} < +\infty$ for some $t > 0$. Let $f : \mathbb{R}^d \to \mathbb{R}$ and let there exist $b_k \in \mathbb{R}$, $k = 1, \ldots, d$, such that $|f(x + ue_k) - f(x)| \leq b_k|u|$ for all $u \in \mathbb{R}$, $x \in \mathbb{R}^d$. Let*

$$h_f(t) = \sup_{x \in \mathbb{R}^d} \sum_{k=1}^d \int_{\mathbb{R}} |f(x + ue_k) - f(x)|^2 \frac{e^{tb_k|u|} - 1}{b_k|u|} \tilde{\nu}_k(\mathrm{d}u), \qquad 0 \leq t < M,$$

*where $M = \sup\{t > 0 : \forall k = 1, \ldots, d, \mathbb{E} e^{tb_k|X_k|} < +\infty\}$. Then,*

$$\mathbb{P}(f(X) - \mathbb{E} f(X) \geq x) \leq e^{-\int_0^x h_f^{-1}(s) \, \mathrm{d}s} \tag{2.1}$$

*for all $0 < x < h_f^{-1}(M^-)$.*

**Proof of Proposition 1 and Corollary 1.** Below, and throughout, by "$f$ Lipschitz with constant $a$" we mean that $|f(x) - f(y)| \leq a\|x - y\|$ for all $x, y \in \mathbb{R}^d$ (the Lipschitz convention stated in [6] also applies). Let us start by recalling the following simple lemma which will be crucial to our approach [5] (or [6] for a sketch of proof). The lemma is the infinitely divisible version of the covariance representation for functions of Gaussian vectors obtained via Gaussian interpolation. Its proof is also obtained via infinitely divisible interpolation and, below, the law of the vector $(U, V)$ is as in the previous proposition.

**Lemma 1.** *Let $X \sim ID(\gamma, 0, \nu)$ be such that $\mathbb{E}\|X\|_2^2 < +\infty$. Let $f, g : \mathbb{R}^d \to \mathbb{R}$ be Lipschitz functions. Then,*

$$\mathbb{E}f(X)g(X) - \mathbb{E}f(X)\mathbb{E}g(X) \tag{2.2}$$
$$= \int_0^1 \mathbb{E}_z\left[\int_{\mathbb{R}^d}(f(U+u) - f(U))(g(V+u) - g(V))\nu(\mathrm{d}u)\right]\mathrm{d}z,$$

*where $\mathbb{E}_z$ is as in Proposition 1.*

We then follow [6]. First, by independence,

$$C = \{t > 0 : \forall k = 1, \ldots, d, \mathbb{E}\mathrm{e}^{tb_k|X_k|} < +\infty\}$$
$$= \left\{t > 0 : \forall k = 1, \ldots, d, \int_{|u|>1} \mathrm{e}^{tb_k|u|}\tilde{\nu}_k(\mathrm{d}u) < +\infty\right\}.$$

Next, we apply the covariance representation (2.2) to $f$ satisfying the above hypotheses and moreover assumed to be bounded and such that $\mathbb{E}f = 0$. Hence,

$$\mathbb{E}f\mathrm{e}^{tf} = \int_0^1 \mathbb{E}_z\left[\mathrm{e}^{tf(V)}\sum_{k=1}^d \int_{\mathbb{R}}(f(U+ue_k) - f(U))(\mathrm{e}^{t(f(V+ue_k)-f(V))} - 1)\tilde{\nu}_k(\mathrm{d}u)\right]\mathrm{d}z$$

$$\leq \int_0^1 \mathbb{E}_z\left[\mathrm{e}^{tf(V)}\sum_{k=1}^d \int_{\mathbb{R}}|f(U+ue_k) - f(U)|\right.$$
$$\left.\times |f(V+ue_k) - f(V)|\frac{\mathrm{e}^{tb_k|u|} - 1}{b_k|u|}\tilde{\nu}_k(\mathrm{d}u)\right]\mathrm{d}z$$

$$\leq \int_0^1 \mathbb{E}_z\left[\mathrm{e}^{tf(V)}\sum_{k=1}^d \int_{\mathbb{R}}\frac{|f(U+ue_k) - f(U)|^2 + |f(V+ue_k) - f(V)|^2}{2}\right.$$
$$\left.\times \left(\frac{\mathrm{e}^{tb_k|u|} - 1}{b_k|u|}\right)\tilde{\nu}_k(\mathrm{d}u)\right]\mathrm{d}z,$$

which gives Proposition 1. For Corollary 1, we continue.

$$\mathbb{E}f\mathrm{e}^{tf} \leq h_f(t)\mathbb{E}[\mathrm{e}^{tf}],$$

where we have used the "marginal property" mentioned above and the fact that $h_f(t)$ is well defined for $0 \leq t < M$. Integrating this last inequality, applied to $f - \mathbb{E}f$, leads to

$$\mathbb{E}e^{t(f-\mathbb{E}f)} \leq e^{\int_0^t h_f(s)\,ds}, \qquad 0 \leq t < M, \tag{2.3}$$

for all bounded $f$ satisfying the hypotheses of the theorem. Fatou's lemma allows us to remove the boundedness assumption in (2.3).

To obtain the tail inequality (2.1), the Bienaymé–Chebyshev inequality gives

$$\mathbb{P}(f(X) - \mathbb{E}f(X) \geq x) \leq \exp\biggl(-\sup_{0<t<M}\biggl(tx - \int_0^t h_f(s)\,ds\biggr)\biggr) = e^{-\int_0^x h_f^{-1}(s)\,ds}$$

by standard arguments (see, e.g., [6]). $\square$

### 2.2. First applications

In general, the above corollary does not provide dimension-free results, even if it slightly improves a result of [6]. However, for particular functions, the above formula can, in fact, be quite efficient. As a consequence of the previous corollary, we present some almost dimension-free results. First, we have the following.

**Theorem 4.** *Let $X$ be as in Theorem 1. Let $\varepsilon > 0$. Then, for all $0 < x < h(M^-)$,*

$$\mathbb{P}(\|X\|_2 \geq (1+\varepsilon)\mathbb{E}\|X\|_2 + x) \leq e^{-\int_0^x h^{-1}(s)\,ds}, \tag{2.4}$$

*where the (dimension-free) function $h$ is given by*

$$h(t) = 8\int_\mathbb{R} |u|(e^{t|u|} - 1)\tilde{\nu}(du) + \frac{2d}{(\varepsilon\mathbb{E}\|X\|_2)^2}\int_\mathbb{R} |u|^3(e^{t|u|} - 1)\tilde{\nu}(du).$$

Theorem 4 still has some weak dimension dependency via the term $\varepsilon\mathbb{E}\|X\|_2$ (the expectation and the median playing the same role up to some constant). In particular, it does not precisely recover Talagrand's result, even for the Euclidean norm. However, the function $h$ itself is dimension-free, in that it can be both upper and lower bounded independently of the dimension $d$ since for $X = (X_1, \ldots, X_d)$, $d\min_{i=1,\ldots,d}(\mathbb{E}|X_i|)^2 \leq (\mathbb{E}\|X\|_2)^2 \leq d\max_{i=1,\ldots,d}\mathbb{E}(X_i^2)$.

The advantage of Theorem 4 is that it does not require any additional assumptions, in contrast to Theorem 2, and that it recovers the $x\log x$-type bound when $\nu$ has bounded support, which is not the case of Theorem 1. We refer the reader to [7], whose results are sometimes superseded by the present paper, for various applications of Theorem 4.

Actually, the mild dimension dependency in Theorem 4 is not much of a problem in the statistical applications we have in mind. However, a statistician would prefer not to have any unnecessary extra assumptions on the variables themselves. Let us explain these comments by means of an example.

Assume that we observe $n$ Poisson processes on $[0,1]$ with intensity $s$ with respect to the Lebesgue measure. We would like to estimate the function $s$. A simple way to do it is to discretize the problem. So, let $d$ be some integer (usually smaller than $\sqrt{n}$ if $s$ is regular enough) and for all $i$, $1 \leq i \leq d$, let $N_i$ be the total number of points that have appeared between $(i-1)/d$ and $i/d$. Therefore, the variables $N_i$ are independent and $N_i$ obeys a Poisson law with parameter $S_i = \int_{(i-1)/d}^{i/d} s(x) n \, dx \simeq (n/d) s(i/d)$. If we want to understand the behavior of the estimator $N = (N_1, \ldots, N_d)$ around $S = (S_1, \ldots, S_d)$, we need to control $\|\epsilon\|_2$, where we write, for all $i$,

$$N_i = S_i + \epsilon_i.$$

We deal here with a regression problem where the noise $\epsilon = (\epsilon_1, \ldots, \epsilon_n)$ has independent components, but these components are not identically distributed. More generally, we would like to encompass the case where the noise is centered with independent infinitely divisible components. Classical regression corresponds to a Gaussian i.i.d. noise and we refer the interested reader to [8] for an extensive study of the link between concentration and estimation of the signal $S$ through model selection methods in that framework.

We would not only like to drop the i.i.d. assumption in Theorem 4 (which is rather easy to do), but we would also like to have an inequality which is valid under a very mild assumption on the noise. The next corollary only assumes that there exists a known bound on the support of the Lévy measures, which is, for instance, the case for the Poisson problem we described above.

**Corollary 2.** *Let $X \sim ID(\gamma, 0, \nu)$ have independent components and be such that $R_k = \inf\{\rho > 0, \tilde{\nu}_k(|x| > \rho) = 0\}$ is finite, with $R = \max_{1 \leq k \leq d} R_k$. Let $\varepsilon > 0$ and let*

$$V_\varepsilon^2 = 8 \max_{1 \leq k \leq d} \int u^2 \tilde{\nu}_k(du) + \frac{2}{(\varepsilon \mathbb{E}\|X\|_2)^2} \sum_{k=1}^d \int u^4 \tilde{\nu}_k(du).$$

*Then, for all $x \geq 0$,*

$$\mathbb{P}(\|X\|_2 \geq (1+\varepsilon)\mathbb{E}\|X\|_2 + x) \leq e^{x/R - (x/R + V_\varepsilon^2/R^2)\log(1 + Rx/V_\varepsilon^2)}.$$

The above result improves upon known concentration inequalities for Poisson processes. One can easily prove that in the framework mentioned above with $X = N - S$, the factor $V_\varepsilon^2$ appearing in Corollary 2 is of the order $8nB/d$, where $B$ is an upper bound on $s$, as soon as $s$ is bounded from below and $d << \sqrt{n}$. So, applying Corollary 2 to $N - S$, we obtain that there exists a constant $C > 0$ such that for all $x \geq 0$,

$$\mathbb{P}(\|N - S\|_2 \geq (1+\varepsilon)\mathbb{E}\|N - S\|_2 + x) \leq e^{-C \min(x^2 d/Bn, x\log(xd/Bn))}.$$

Applying instead Theorem 4 of [10] to our problem, we see that there exists a constant $C' > 0$ such that for all $x \geq 0$,

$$\mathbb{P}(\|N - S\|_2 \geq (1+\varepsilon)\mathbb{E}\|N - S\|_2 + x) \leq e^{-C' \min(x^2 d/Bn, x)}.$$

Thus, even in this simplest case (of a Poisson process), we see that Corollary 2 (optimally) improves known results by a logarithmic factor.

However, if one is interested in recovering the function $s$ from the observations of the $n$ Poisson processes and if $s$ is not very smooth and varies greatly on a very small interval, looking at a regular partition might be rather useless. It might, instead, be much more fruitful to look at the function $s$ discretized on very small intervals ($d = n$) and then to look at the projection of the signal $S$ on a space $\mathcal{S}$ which is generated by, say, a few Haar wavelets. Now, if the signal $S$ is sparse on the Haar basis, it means that one can find (and use) a space $\mathcal{S}$, with a dimension much smaller than $n$, to provide a good approximation for $S$. However, one now needs to understand the behavior of $\|\Pi_{\mathcal{S}}(N - S)\|_2$, where $\Pi_{\mathcal{S}}$ is the orthogonal projection on $\mathcal{S}$. As before, we want to deal with more general infinitely divisible noise than simply centered Poisson variables. The two following corollaries provide such results.

**Corollary 3.** *Let $X \sim ID(\gamma, 0, \nu)$ have independent components. Let $\mathcal{S}$ be a subspace of $\mathbb{R}^d$ and $\Pi_{\mathcal{S}}$ the orthogonal projection on $\mathcal{S}$. Let $M = \sup\{t > 0 : \forall k = 1, \ldots, d, \mathbb{E}e^{t|X_k|} < +\infty\}$. Let $E > 0$. Then, for all $0 < x < h(M^-)$,*

$$\mathbb{P}(\|\Pi_{\mathcal{S}}(X)\|_2 \geq \mathbb{E}\|\Pi_{\mathcal{S}}(X)\|_2 + E + x) \leq e^{-\int_0^x h^{-1}(s)\mathrm{d}s} \qquad (2.5)$$

*and*

$$\mathbb{P}(\|\Pi_{\mathcal{S}}X\|_2 \leq \mathbb{E}\|\Pi_{\mathcal{S}}(X)\|_2 - E - x) \leq e^{-\int_0^x h^{-1}(s)\mathrm{d}s}, \qquad (2.6)$$

*where the function $h$ is given by*

$$h(t) = 8 \max_{1 \leq k \leq d} \int_{\mathbb{R}} |u|(e^{t|u|} - 1)\tilde{\nu}_k(\mathrm{d}u) + \frac{2}{E^2} \sum_{k=1}^d \|\Pi_{\mathcal{S}}(e_k)\|_2^4 \int_{\mathbb{R}} |u|^3(e^{t|u|} - 1)\tilde{\nu}_k(\mathrm{d}u)$$

*for $0 \leq t < M$.*

The next version, which assumes i.i.d. coordinates, in contrast to the above, can sometimes be easier to use. For the statistician, the i.i.d. case appears when dealing with a regression problem where the noise does not depend on the signal itself.

**Corollary 4.** *Let $X$ be as in Theorem 1. Let $\mathcal{S}$ be a subspace of $\mathbb{R}^d$ and let $\Pi_{\mathcal{S}}$ be the orthogonal projection on $\mathcal{S}$. Let $\varepsilon > 0$. Then, for all $0 < x < h(M^-)$,*

$$\mathbb{P}(\|\Pi_{\mathcal{S}}(X)\|_2 \geq (1 + \varepsilon)\sqrt{\mathbb{E}\|\Pi_{\mathcal{S}}(X)\|_2^2} + x) \leq e^{-\int_0^x h^{-1}(s)\mathrm{d}s} \qquad (2.7)$$

*and*

$$\mathbb{P}(\|\Pi_{\mathcal{S}}(X)\|_2 \leq \mathbb{E}\|\Pi_{\mathcal{S}}(X)\|_2 - \varepsilon\sqrt{\mathbb{E}\|\Pi_{\mathcal{S}}(X)\|_2^2} - x) \leq e^{-\int_0^x h^{-1}(s)\mathrm{d}s}, \qquad (2.8)$$

*where the (dimension-free) function h is given by*

$$h(t) = 8\int_{\mathbb{R}} |u|(e^{t|u|} - 1)\tilde{\nu}(\mathrm{d}u) + \frac{2}{\varepsilon^2 \mathbb{E}X_1^2} \int_{\mathbb{R}} |u|^3 (e^{t|u|} - 1)\tilde{\nu}(\mathrm{d}u)$$

*for $0 \leq t < M$.*

Finally, in the density framework [4], it is known that the Euclidean norm does not suffice to assess the performance of one estimator: if the density $s$ belongs to a Sobolev space $H^\alpha$, linear estimators cannot achieve the optimal rate of convergence for the $L_p$-norm if $p > 2$, whereas they can for the $L_2$-norm. This corresponds to sparse signals $S$ that can be approximated by their projection on a subspace with small dimension with respect to $n$ but for the $L_p$-norm. So, in that context, it is necessary to work with general $L_p$-norms and concentration results for $\ell_p$-norms are more relevant than for the $\ell_2$-norm. Such results are presented in the following corollary with a weak dimension dependency and in Theorem 5 without any dimension dependency and with as few assumptions as possible on the noise (see Section 4). These results are given in their i.i.d. version for the sake of simplicity, but one can easily see their non-i.i.d. version from the given proofs.

**Corollary 5.** *Let $X$ be as in Theorem 1. Let $p \geq 2$ and $\varepsilon > 0$. Then, for all $0 < x < h(M^-)$,*

$$\mathbb{P}(\|X\|_p \geq (1+\varepsilon)\mathbb{E}(\|X\|_p) + x) \leq e^{-\int_0^x h^{-1}(s)\mathrm{d}s} \tag{2.9}$$

*and*

$$\mathbb{P}(\|X\|_p \leq (1-\varepsilon)\mathbb{E}(\|X\|_p) - x) \leq e^{-\int_0^x h^{-1}(s)\mathrm{d}s}, \tag{2.10}$$

*where the function $h$ is given by*

$$h(t) = p^2 \int_{\mathbb{R}} \left(1 + \frac{|u|d^{1/(2p-2)}}{\varepsilon \mathbb{E}(\|X\|_p)}\right)^{2p-2} |u|(e^{t|u|} - 1)\tilde{\nu}(\mathrm{d}u)$$

*for $0 \leq t < M$.*

### 2.3. Proofs

Let us proceed to the proof of the results of the previous subsection. We begin with Corollary 3, the other proofs being easier.

**Proof of Corollary 3.** We apply Corollary 1 to $f(x) = (\|\Pi_S(x)\|_2 - E)^+$. First, it is easily verified that for each $k$, $|f(x + ue_k) - f(x)| \leq |\|\Pi_S(x + ue_k)\|_2 - \|\Pi_S(x)\|_2 |\mathbf{1}_{A_k}$, where $A_k = \{\|\Pi_S(x + ue_k)\|_2 \geq E \text{ or } \|\Pi_S(x)\|_2 \geq E\}$. We then have

$$|f(x + ue_k) - f(x)| \leq \frac{|2\langle u\Pi_S(e_k), \Pi_S(x)\rangle + u^2\|\Pi_S(e_k)\|_2^2|\mathbf{1}_{A_k}}{\|\Pi_S(x + ue_k)\|_2 + \|\Pi_S(x)\|_2}$$

$$
\leq \frac{2|u||\langle \Pi_S(e_k), \Pi_S(x)\rangle|}{\|\Pi_S(x)\|_2} + \frac{u^2 \|\Pi_S(e_k)\|_2^2}{E}. \tag{2.11}
$$

Moreover, since $|f(x+ue_k) - f(x)| \leq |u|$, we have

$$
\sum_{k=1}^{d} \int_{\mathbb{R}} |f(x+ue_k) - f(x)|^2 \frac{e^{tb_k|u|} - 1}{b_k|u|} \tilde{\nu}_k(\mathrm{d}u)
$$

$$
\leq \sum_{k=1}^{d} \int_{\mathbb{R}} \left( 8u^2 \frac{|\langle \Pi_S(e_k), \Pi_S(x)\rangle|^2}{\|\Pi_S(x)\|_2^2} + \frac{2u^4 \|\Pi_S(e_k)\|_2^4}{E^2} \right) \left( \frac{e^{t|u|} - 1}{|u|} \right) \tilde{\nu}_k(\mathrm{d}u)
$$

$$
\leq \sum_{k=1}^{d} \int_{\mathbb{R}} \left( 8u^2 \frac{|\langle e_k, \Pi_S(x)\rangle|^2}{\|\Pi_S(x)\|_2^2} + \frac{2u^4 \|\Pi_S(e_k)\|_2^4}{E^2} \right) \left( \frac{e^{t|u|} - 1}{|u|} \right) \tilde{\nu}_k(\mathrm{d}u).
$$

Hence, $h_f \leq h$. To complete the proof of (2.5), note that $\|\Pi_S(X)\|_2 - E \leq (\|\Pi_S(X)\|_2 - E)^+$ and that $\mathbb{E}(\|\Pi_S(X)\|_2 - E)^+ \leq \mathbb{E}\|\Pi_S(X)\|_2$. To get the lower bound (2.6), just proceed as above, but with the function $f(x) = -(\|\Pi_S(x)\|_2 - E)^+$ and note that $(\|\Pi_S(X)\|_2 - E)^+ \leq \|\Pi_S(X)\|_2$ and that $\mathbb{E}\|\Pi_S(X)\|_2 - E \leq \mathbb{E}(\|\Pi_S(X)\|_2 - E)^+$. □

**Proof of Theorem 4.** We apply Corollary 3 with $S = \mathbb{R}^d$ and $E = \varepsilon \mathbb{E}\|X\|_2$. □

**Proof of Corollary 2.** It is sufficient to note that, proceeding as in Corollary 3, $h(t) \leq h_0(t) = V_\varepsilon^2(e^{xR} - 1/R)$ and that $M = +\infty$. It remains to integrate the reciprocal of $h_0$. □

**Proof of Corollary 4.** Again applying Corollary 3, let us take $E = \varepsilon \sqrt{\mathbb{E}(\|\Pi_S(X)\|_2^2)}$. Then, note that in the centered i.i.d. case,

$$
\mathbb{E}[\|\Pi_S(X)\|_2^2] = \mathbb{E}\left[ \sum_{l=1}^{d} \left( \sum_{k=1}^{d} X_k \langle \Pi_S(e_k), e_l \rangle \right)^2 \right]
$$

$$
= \sum_{l=1}^{d} \sum_{k=1}^{d} \mathbb{E}[X_k^2] \langle \Pi_S(e_k), e_l \rangle^2
$$

$$
= \mathbb{E}[X_1^2] \sum_{k=1}^{d} \|\Pi_S(e_k)\|_2^2
$$

$$
\geq \mathbb{E}[X_1^2] \sum_{k=1}^{d} \|\Pi_S(e_k)\|_2^4
$$

since $\|\Pi_S(e_k)\|_2 \leq \|e_k\|_2 = 1$. □

**Proof of Corollary 5.** We apply Corollary 1 to $f(x) = (\|x\|_p - \varepsilon \mathbb{E}(\|X\|_p))^+$ to get the first result. For the second one, it is sufficient to do the same with $-(\|x\|_p - \varepsilon \mathbb{E}(\|X\|_p))^+$.

The end of the proof is as in Corollary 3. Let $x$ be in $\mathbb{R}^d$. First, it is easily verified that for each $k$, $|f(x+ue_k) - f(x)| \le |\|x+ue_k\|_p - \|x\|_p| \mathbf{1}_{A_k}$, where $A_k = \{\|x+ue_k\|_p \ge \varepsilon \mathbb{E}(\|X\|_p)$ or $\|x\|_p \ge \varepsilon \mathbb{E}(\|X\|_p)\}$. Since

$$\forall a, b \ge 0 \qquad |a-b| \le \frac{|a^p - b^p|}{\sup(a,b)^{p-1}}, \tag{2.12}$$

we have

$$|f(x+ue_k) - f(x)| \le \frac{||x_k+u|^p - |x_k|^p|}{\sup(\|x\|_p, \|x+ue_k\|_p)^{p-1}} \mathbf{1}_{A_k}. \tag{2.13}$$

But since $x \mapsto x^p$ is convex, we have

$$||x_k + u|^p - |x_k|^p| \le |(|x_k| + |u|)^p - |x_k|^p|.$$

Combining this with the fact that

$$\forall y \ge 0 \qquad (1+y)^p - 1 \le py(1+y)^{p-1}$$

implies that

$$\begin{aligned}|f(x+ue_k) - f(x)| &\le \frac{p|u|(|x_k|+|u|)^{p-1}}{\sup(\|x\|_p, \|x+ue_k\|_p)^{p-1}} \mathbf{1}_{A_k} \\ &\le \frac{p|u|(|x_k|+|u|)^{p-1}}{\sup(\|x\|_p, \varepsilon\mathbb{E}(\|X\|_p))^{p-1}}.\end{aligned} \tag{2.14}$$

Moreover, since $|f(x+ue_k) - f(x)| \le |u|$, we have

$$\begin{aligned}\sum_{k=1}^d \int_{\mathbb{R}} &|f(x+ue_k) - f(x)|^2 \frac{e^{tb_k|u|} - 1}{b_k|u|} \tilde{\nu}_k(\mathrm{d}u) \\ &\le p^2 \int_{\mathbb{R}} \frac{\||x| + |u|I\|_{2p-2}^{2p-2}}{\sup(\|x\|_p, \varepsilon\mathbb{E}(\|X\|_p))^{2p-2}} |u|(e^{t|u|} - 1)\tilde{\nu}(\mathrm{d}u) \\ &\le p^2 \int_{\mathbb{R}} \left(\frac{\|x\|_{2p-2} + |u|\|I\|_{2p-2}}{\sup(\|x\|_p, \varepsilon\mathbb{E}(\|X\|_p))}\right)^{2p-2} |u|(e^{t|u|} - 1)\tilde{\nu}(\mathrm{d}u),\end{aligned} \tag{2.15}$$

where $I = (1,\ldots,1) \in \mathbb{R}^d$ and $|x| = (|x_1|,\ldots,|x_d|)$. Since $p \ge 2$, $2p-2 \ge p$ and $\|x\|_{2p-2} \le \|x\|_p$, which implies that

$$\begin{aligned}\sum_{k=1}^d \int_{\mathbb{R}} &|f(x+ue_k) - f(x)|^2 \frac{e^{tb_k|u|} - 1}{b_k|u|} \tilde{\nu}_k(\mathrm{d}u) \\ &\le p^2 \int_{\mathbb{R}} \left(1 + |u|\frac{d^{1/(2p-2)}}{\varepsilon\mathbb{E}(\|X\|_p)}\right)^{2p-2} |u|(e^{t|u|} - 1)\tilde{\nu}(\mathrm{d}u).\end{aligned}$$

Here, again, the upper bound is dimension-free since $\mathbb{E}(\|X\|_p) \geq \mathbb{E}(|X_1|)d^{1/p}$. □

## 3. Proof of Theorem 1

Another method to decouple $U$ and $V$ in Proposition 1 is to use the following inequality, which is a particular instance of Young's inequality (for the pairs of conjugate functions $ce^x$ and $y\log(y/c) - y$, with the optimal $c$) and has already been used in [8].

**Lemma 2.** *Let $\lambda > 0$ and let $X$ and $Y$ be random variables for which all the expectations below exist. Then,*

$$\mathbb{E}[X e^{\lambda Y}] \leq \mathbb{E}[Y e^{\lambda Y}] + \frac{\log \mathbb{E}[e^{\lambda X}]}{\lambda} \mathbb{E}[e^{\lambda Y}] - \frac{\log \mathbb{E}[e^{\lambda Y}]}{\lambda} \mathbb{E}[e^{\lambda Y}]. \tag{3.1}$$

**Proof.** Indeed, if

$$dQ = \frac{e^{\lambda Y}}{\mathbb{E}[e^{\lambda Y}]} d\mathbb{P},$$

then, by Jensen's inequality,

$$\lambda \mathbb{E}_Q(X - Y) \leq \log \mathbb{E}_Q(e^{\lambda(X-Y)}). \qquad \square$$

With the help of the previous lemma, the following holds.

**Corollary 6.** *Let $X = (X_1, \ldots, X_d) \sim ID(\gamma, 0, \nu)$ have i.i.d. components and be such that $\mathbb{E}e^{t\|X\|_2} < +\infty$ for some $t > 0$. Let $f: \mathbb{R}^d \to \mathbb{R}$ be such that $\mathbb{E}f(X) = 0$ and let there exist $b \in \mathbb{R}$ such that for all $k$, $|f(x + ue_k) - f(x)| \leq b|u|$ for all $u \in \mathbb{R}$, $x \in \mathbb{R}^d$. Assume, moreover, that for all $u \in \mathbb{R}$, there exists a function $C_u$ such that*

$$\sum_{k=1}^{d} \int_{\mathbb{R}} |f(X + ue_k) - f(X)|^2 \tilde{\nu}_k(du) \leq u^2 C_u(X)$$

*and such that $\mathbb{E}[e^{\lambda(u,t)C_u(X)}] < \infty$ for $\lambda(u,t) > 0$. Then, for all $t$ for which all the quantities below are well defined, we have*

$$(1 - h(t))\mathbb{E}[f e^{tf}] \leq g(t)\mathbb{E}[e^{tf}], \tag{3.2}$$

*where*

$$h(t) = \int_{\mathbb{R}} \frac{t}{\lambda(u,t)} |u| \frac{e^{tb|u|} - 1}{b} \tilde{\nu}(du)$$

*and*

$$g(t) = \int_{\mathbb{R}} \frac{\ln(\phi(u,t))}{\lambda(u,t)} |u| \frac{e^{tb|u|} - 1}{b} \tilde{\nu}(du),$$

and where $\phi(u,t) = \mathbb{E}[e^{\lambda(u,t)C_u(X)}]$.

**Proof.** Applying Proposition 1 to $f$, the above assumptions entail that

$$\mathbb{E}fe^{tf} \leq \int_0^1 \mathbb{E}_z \left[ \int_{\mathbb{R}} \frac{C_u(U) + C_u(V)}{2} e^{tf(V)} \left( |u| \frac{e^{tb|u|} - 1}{b} \right) \tilde{\nu}(du) \right] dz.$$

Next, apply Lemma 2 to $\lambda(u,t)Y = tf(V)$ and to $X = C_u(U)$ or $X = C_u(V)$. Since $Y$ has zero mean, one can ignore the last term in (3.1), and this leads to

$$\mathbb{E}fe^{tf} \leq \int_{\mathbb{R}} \left[ \mathbb{E}\left( \frac{t}{\lambda(u,t)} fe^{tf} \right) + \frac{\ln(\phi(u,t))}{\lambda(u,t)} \mathbb{E}(e^{tf}) \right] |u| \frac{e^{tb|u|} - 1}{b} \tilde{\nu}(du),$$

which concludes the proof. □

(3.2) is non-trivial only when $h(t) < 1$. One of its applications is to prove Theorem 1, with the help of our next two results.

**Lemma 3.** *For $\alpha > 0$, let $\ell_\alpha = -\ln \mathbb{E}[e^{-\alpha X_1^2}]$. Then, for all $\lambda > 0, v > 0$ such that, $\ell_\alpha \geq \lambda/v$,*

$$\mathbb{E}\left( \exp\left( \frac{\lambda d}{\|X\|_2^2 + v} \right) \right) \leq 1 + \exp\left( \frac{\alpha \lambda}{\ell_\alpha - \lambda/v} \right).$$

**Proof.** Let $\varepsilon > 0$, to be chosen later. Let $a = \exp(\frac{\lambda d}{\varepsilon d \mathbb{E}(X_1^2) + v})$ and let $b = \exp(\frac{\lambda d}{v})$. Then,

$$\mathbb{E}\left( \exp\left( \frac{\lambda d}{\|X\|_2^2 + v} \right) \right) = \int_0^b \mathbb{P}\left( \exp\left( \frac{\lambda d}{\|X\|_2^2 + v} \right) \geq t \right) dt$$

$$\leq a + \int_a^b \mathbb{P}\left( -\|X\|_2^2 \geq v - \frac{\lambda d}{\ln t} \right) dt$$

$$\leq a + \int_a^b \mathbb{E}[e^{-\alpha X_1^2}]^d e^{-\alpha v + \alpha \lambda d / \ln t} dt$$

$$\leq a + e^{-d\ell_\alpha + \alpha \varepsilon d \mathbb{E}(X_1^2)}(b - a)$$

$$\leq a + e^{d(\varepsilon \alpha \mathbb{E}(X_1^2) + (\lambda/v) - \ell_\alpha)}.$$

Taking $\varepsilon$ such that $\varepsilon \alpha \mathbb{E}(X_1^2) + (\lambda/v) - \ell_\alpha = 0$ leads to

$$\mathbb{E}\left( \exp\left( \frac{\lambda d}{\|X\|_2^2 + v} \right) \right) \leq a + 1 \leq 1 + \exp\left( \frac{\alpha \lambda}{\ell_\alpha - \lambda/v} \right). \qquad \square$$

**Lemma 4.** *There exists positive constants $c_1, c_2, c_3$ such that for all $x \in \mathbb{R}^d$ and $u \in \mathbb{R}$,*

$$\sum_{k=1}^{d} |\|x + ue_k\|_2 - \|x\|_2|^2 \leq u^2\left(c_1 + \frac{c_2\, du^2}{\|x\|_2 + c_3 u^2}\right).$$

**Proof.** The proof is similar to an argument used in the proof of Corollary 3. We have

$$\sum_{k=1}^{d} |\|x + ue_k\|_2 - \|x\|_2|^2 \leq \sum_{k=1}^{d}\left(\frac{2ux_k + u^2}{\|x + ue_k\|_2 + \|x\|_2}\right)^2.$$

But, for all $\varepsilon > 0$,

$$\|x + ue_k\|_2^2 = \sum_{j \neq k} x_j^2 + (x_k + u)^2 = \|x\|_2^2 + 2ux_k + u^2 \geq \|x\|_2^2 - \varepsilon x_k^2 + (1 - \varepsilon^{-1})u^2.$$

Therefore,

$$(\|x + ue_k\|_2 + \|x\|_2)^2 \geq \|x + ue_k\|_2^2 + \|x\|_2^2 \geq (2 - \varepsilon)\|x\|_2^2 + (1 - \varepsilon^{-1})u^2.$$

Taking $\varepsilon = 3/2$, completes the proof. $\square$

With the help of the previous lemmas, we now get the following.

**Proof of Theorem 1.** We want to apply Corollary 6 to $f(X) = \|X\|_2 - \mathbb{E}\|X\|_2$. With the notation of Lemma 4,

$$C_u(X) = c_1 + \frac{c_2\, du^2}{\|X\|_2 + c_3 u^2}$$

works. We must then compute $\ln(\phi(u,t))$. But,

$$\ln(\phi(u,t)) = c_1 \lambda(u,t) + \ln\left(\mathbb{E}\left[\exp\left(\frac{\lambda(u,t)c_2 u^2 d}{\|X\|_2^2 + c_3 u^2}\right)\right]\right)$$

and so from Lemma 3, it follows that

$$\ln(\phi(u,t)) \leq c_1 \lambda(u,t) + \ln\left(1 + \exp\left(\frac{\alpha \lambda(u,t)c_2 u^2}{\ell_\alpha - \lambda(u,t)c_2/c_3}\right)\right)$$

for all $\alpha$ such that $\ell_\alpha > \lambda(u,t)c_2/c_3$.
Taking $\alpha = 1$ and $\lambda(u,t) = c_3 l/(2c_2)$ gives

$$\ln(\phi(u,t)) \leq c_1 \lambda(u,t) + \ln 2 + c_3 u^2,$$

which leads to the result by standard arguments. $\square$

## 4. Proof of Theorem 2

We state and prove in this section a generalization of Theorem 2. First, recall that the vector $X$ can be viewed as the value at time 1 of a Lévy process $(X_z, z \geq 0)$. For every $z \in [0,1]$, decompose $X = Y_z + Z_z$, where $Y_z = X_z$ and $Z_z = X - X_z$, so that $Y_z, Z_z$ are independent. If $1 \leq k \leq d$, we write $(Y_k)_z, (Z_k)_z$ for the $k$th coordinate of $Y_z, Z_z$. Furthermore, to simplify the notation, we denote

$$Y_t^+ = (Y_1)_t \mathbf{1}_{(Y_1)_t \geq 0} \tag{4.1}$$

and likewise define $Y_t^-$, $Z_t^+$, $Z_t^-$.

**Theorem 5.** *Let $X$ be as in Theorem 1 and let $2 \leq p < \infty$. Then, for all $0 < x < h_p(M^-)$,*

$$\mathbb{P}(\|X\|_p - \mathbb{E}\|X\|_p \geq x) \leq \exp\left(-\int_0^x h_p^{-1}(s)\,\mathrm{d}s\right), \tag{4.2}$$

*where the (dimension-free) function $h_p$ is given as follows:*

- *if $X$ has almost surely non-negative coordinates,*

$$h_p(t) = p^2 \int_{\mathbb{R}_+} \left[\left(\frac{1}{2^{1/p}} + \frac{u}{m_p^{1/p}}\right)^{2p-2} + 4\left(1 + \frac{1}{2^{1/p}}\right)^{2p-2} \frac{m_{2p}}{m_p^2}\right] u(e^{tu} - 1)\tilde{\nu}(\mathrm{d}u), \tag{4.3}$$

 *where the moments $m_{2p}, m_p$ are defined by $m_q = \mathbb{E}[X_1^q]$ for $q = p, 2p$;*
- *in the general case,*

$$h_p(t) = p^2 \int_{\mathbb{R}} \left[\left(1 + \frac{2^{1/p}|u|}{\underline{m_p}^{1/p}}\right)^{2p-2} + 2^{2p}\frac{\overline{m_{2p}}}{\underline{m_p}^2}\right] |u|(e^{t|u|} - 1)\tilde{\nu}(\mathrm{d}u), \tag{4.4}$$

*where the modified moments $\overline{m_{2p}}, \underline{m_p}$ are defined, using the notation given in (4.1), by*

$$\underline{m_p} = \inf_{z \in [0,1]}[\inf\{\mathbb{E}[|Y_z^+ + Z_z^+|^p], \mathbb{E}[|Y_z^- + Z_z^-|^p]\}]$$

*and*

$$\overline{m_{2p}} = \sup_{z \in [0,1]}[\sup\{\mathbb{E}[|Y_z^+ + Z_z^+|^{2p}], \mathbb{E}[|Y_z^- + Z_z^-|^{2p}]\}].$$

*When $1 \leq p < 2$, an inequality similar to (4.2) holds, where $h_p$ is now replaced by the following function $h_{p,d}$ (which is no longer dimension-free):*

$$h_{p,d}(t) = p^2 \int_{\mathbb{R}} \left[d^{2/p-1}\left(1 + \frac{2^{1/p}|u|}{\underline{m_p}^{1/p}}\right)^{2p-2} + 2^{2p}\frac{\overline{m_{2p}}}{\underline{m_p}^2}\right] |u|(e^{t|u|} - 1)\tilde{\nu}(\mathrm{d}u).$$

Observe that Theorem 2 is a particular case of this more general result. The obvious drawback of Theorem 5 is that, except in the cases considered in Theorem 2, we do not have a precise control of the quantities $m_p$, $\overline{m_{2p}}$. In particular, if $m_p = 0$, then $h_p = \infty$ and we get a trivial bound. However, it should be clear that when $x$ does not have almost surely positive coordinates, the quantity

$$\inf_{z \in [\varepsilon, 1]} [\inf\{\mathbb{E}[|Y_z^+ + Z_z^+|^p], \mathbb{E}[|Y_z^- + Z_z^-|^p]\}]$$

is positive for every $\varepsilon > 0$. So, the only case when $m_p$ might be zero is the case when $(X_1)_t$, the first coordinate of the Lévy process $X_1$, taken at time $t$, has a probability tending to 1 or to 0 to be positive as $t \to 0$. This kind of behavior does exist, but in most 'natural' examples, this is not the case and then Theorem 5 does give a non-trivial bound, although the expression of this bound is not always easy to handle.

**Proof of Theorem 5.** The proof can be divided into three steps. Define the function $f(V) = \|V\|_p - \mathbb{E}\|X\|_p$.

**Step 1**: We claim that

$$|f(X + ue_k) - f(X)| \leq p|u| \left( \frac{|X_k| + |u|}{(\|X\|_p^p + |u|^p)^{1/p}} \right)^{p-1}. \tag{4.5}$$

**Step 2**: Using the notation of Proposition 1, we have

$$\mathcal{A} := \mathbb{E}_z \left[ \sum_{k=1}^d |f(V + ue_k) - f(V)|^2 e^{tf(V)} \right] \leq p^2 |u|^2 F(p, d, u) \mathbb{E}[e^{tf(V)}],$$

$$\mathcal{B} := \mathbb{E}_z \left[ \sum_{k=1}^d |f(U + ue_k) - f(U)|^2 e^{tf(V)} \right] \leq p^2 |u|^2 F(p, d, u) \mathbb{E}[e^{tf(V)}]$$

for some function $F(p, d, u)$ that will be made explicit in the course of the proof.

**Step 3**: For $h_p$ and $h_{p,d}$ as in Theorem 5, we have, if $p \geq 2$,

$$\mathbb{E}(fe^{tf}) \leq h_p(t)\mathbb{E}(e^{tf})$$

and if $1 \leq p < 2$,

$$\mathbb{E}(fe^{tf}) \leq h_{p,d}(t)\mathbb{E}(e^{tf}).$$

Integrating the inequalities of this last step leads to Theorem 5.

**Proof of Step 1.** First, for all reals $a, b \geq 0$,

$$\max(a,b)^{1/p} - \min(a,b)^{1/p} = \left( \frac{a+b}{2} + \frac{|a-b|}{2} \right)^{1/p} - \left( \frac{a+b}{2} - \frac{|a-b|}{2} \right)^{1/p}$$

$$= \left( \frac{a+b}{2} \right)^{1/p} \left[ \left( 1 + \frac{|a-b|}{a+b} \right)^{1/p} - \left( 1 - \frac{|a-b|}{a+b} \right)^{1/p} \right].$$

Concentration for norms of infinitely divisible vectors     943Since the function $x \mapsto (1+x)^{1/p} - (1-x)^{1/p}$ is convex on $[0,1]$, is zero at zero and $2^{1/p}$ at 1, we get

$$\max(a,b)^{1/p} - \min(a,b)^{1/p} \leq \left(\frac{a+b}{2}\right)^{1/p}\left(2^{1/p}\frac{|a-b|}{a+b}\right) = \frac{|a-b|}{(a+b)^{(p-1)/p}}. \tag{4.6}$$

Next, we want to apply (4.6) to $a = \|X\|_p^p$, $b = \|X + ue_k\|_p^p$. Put $c_k = u/X_k$ and distinguish between three cases. First, if $c_k \geq 0$, then

$$||X_k + u|^p - |X_k|^p| \leq p|u|(|X_k| + |u|)^{p-1}$$

and

$$|X_k + u|^p + |X_k|^p \geq 2|X_k|^p + |u|^p.$$

Second, if $-2 \leq c_k \leq 0$, then

$$||X_k + u|^p - |X_k|^p| \leq p|u||X_k|^{p-1}$$

and we can check that for every $A \geq 0$,

$$\frac{|X_k|}{(A + |X_k + u|^p + |X_k|^p)^{1/p}} \leq \frac{(1 + |c_k|)|X_k|}{(A + |X_k|^p + |c_k|^p|X_k|^p)^{1/p}}.$$

Third, if $c_k \leq -2$, then

$$||X_k + u|^p - |X_k|^p| \leq p|u||X_k + u|^{p-1} = p|u|[|1 + c_k||X_k|]^{p-1}$$

and we can check that for every $A \geq 0$,

$$\frac{|1 + c_k||X_k|}{(A + |X_k + u|^p + |X_k|^p)^{1/p}} \leq \frac{(1 + |c_k|)|X_k|}{(A + |X_k|^p + |c_k|^p|X_k|^p)^{1/p}}.$$

Combining the inequalities in these three cases with (4.6), we obtain (4.5).

**Proof of Step 2.** We will use the following well-known result.

**Lemma 5.** *Let $T$ be a random vector in $\mathbb{R}^d$ with i.i.d. components and let $A, B : \mathbb{R}^d \to \mathbb{R}$ be two functions. For every $i \leq d$ and every $x \in \mathbb{R}^{d-1}$, define the functions $A_{x,i}, B_{x,i} : \mathbb{R} \to \mathbb{R}$ via*

$$A_{x,i}(t) = A(x_1, \ldots, x_{i-1}, t, x_{i+1}, \ldots, x_d),$$
$$B_{x,i}(t) = (x_1, \ldots, x_{i-1}, t, x_{i+1}, \ldots, x_d).$$

*Assume that for every $i \leq d$ and every $x \in \mathbb{R}^{d-1}$, one of the two functions $A_{x,i}, B_{x,i}$ is non-decreasing and the other one is non-increasing. Then*

$$\mathbb{E}[A(T)B(T)] \leq \mathbb{E}[A(T)]\mathbb{E}[B(T)]. \tag{4.7}$$

**Proof.** When $d=1$, the proof is obtained by writing

$$\mathbb{E}[A(T)B(T)] - \mathbb{E}[A(T)]\mathbb{E}[B(T)] = \tfrac{1}{2}\mathbb{E}[(A(T)-A(T'))(B(T)-B(T'))],$$

where $T'$ is an independent copy of $T$, and using the monotonicity assumptions. When $d \geq 1$, we proceed by induction. $\square$

Let us first bound $\mathcal{A}$. Summing (4.5) over $k$, the triangle inequality for $\|\cdot\|_{2p-2}$ gives

$$\mathcal{A} \leq p^2 |u|^2 \mathbb{E}_z\left[\left(\frac{\|V\|_{2p-2} + |u|d^{1/(2p-2)}}{(\|V\|_p^p + |u|^p)^{1/p}}\right)^{2p-2} \mathrm{e}^{tf(V)}\right].$$

If $2p-2 \geq p$, $\|V\|_{2p-2} \leq \|V\|_p$ and, otherwise, $\|V\|_{2p-2} \leq \|V\|_p d^{1/(2p-2)-1/p}$. Consequently,

$$\mathcal{A} \leq p^2 |u|^2 \mathbb{E}_z\left[Q_p\left(\frac{(\|V\|_p^p + |u|^p)^{1/p}}{|u|d^{1/(2p-2)}}\right)\mathrm{e}^{tf(V)}\right],$$

where the polynomials $Q_p$ are defined as follows:

- if $p \geq 2$, $Q_p : x \mapsto (1 + 1/x)^{2p-2}$;
- if $p < 2$, $Q_p : x \mapsto (d^{1/(2p-2)-1/p} + 1/x)^{2p-2}$.

We can apply (4.7) to $T = \|V\|_p$ since $Q_p(\frac{(\|V\|_p^p+|u|^p)^{1/p}}{|u|d^{1/(2p-2)}})$ is decreasing in $\|V\|_p$, while $\mathrm{e}^{tf(V)}$ is increasing. This gives

$$\mathcal{A} \leq p^2 |u|^2 \mathbb{E}\left[Q_p\left(\frac{(\|V\|_p^p + |u|^p)^{1/p}}{|u|d^{1/(2p-2)}}\right)\right]\mathbb{E}[\mathrm{e}^{tf(V)}].$$

Next, let us bound $\mathcal{B}$. Summing (4.5) over $k$ gives $\mathcal{B} \leq p^2 |u|^2 \mathcal{D}$ with

$$\mathcal{D} = \mathbb{E}_z\left[Q_p\left(\frac{(\|U\|_p^p + |u|^p)^{1/p}}{|u|d^{1/(2p-2)}}\right)\mathrm{e}^{tf(V)}\right].$$

Note that the distribution $\mathbb{P}_z$ of $(U,V)$ is such that for all $z \in [0,1]$, $U = Y_z' + Z_z$, while $V = Y_z + Z_z$, where $Y_z, Y_z'$ and $Z_z$ are independent and $Y_z$ and $Y_z'$ are identically distributed. Hence,

$$\mathcal{D} = \mathbb{E}_z\left[Q_p\left(\frac{(\|Y_z' + Z_z\|_p^p + |u|^p)^{1/p}}{|u|d^{1/(2p-2)}}\right)\mathrm{e}^{tf(Y_z + Z_z)}\right].$$

Let us introduce the function $M_p$: for all real $u, y, y', z$ in $\mathbb{R}^d$,

$$M_p(|u|, y, y', z) = \left[|u|^p + \sum_i \mathbf{1}_{sgn(y_i) = sgn(y_i') = sgn(z_i)}|z_i + y_i'|^p\right]^{1/p}.$$

Then, $M_p(|u|, y, y', z) \leq (\|y' + z\|_p^p + |u|^p)^{1/p}$ and since $Q_p$ is decreasing,

$$\mathcal{D} \leq \mathbb{E}_z\left[Q_p\left(\frac{M_p(|u|, Y_z, Y'_z, Z_z)}{|u|d^{1/(2p-2)}}\right)e^{tf(Y_z+Z_z)}\right]. \tag{4.8}$$

From now on, we split $\mathbb{E}_z$ into $\mathbb{E}_{Y_z,Y'_z}\mathbb{E}_{Z_z}$, meaning that we first integrate according to the law of $Z_z$. Let us fix $Y_z$ and $Y'_z$. Suppose, first, that $(Z_k)_z$, $(Y_k)_z$ and $(Y'_k)_z$ have the same sign. Then, one of the two functions $Q_p(\frac{M_p(|u|,Y_z,Y'_z,Z_z)}{|u|d^{1/(2p-2)}})$, $e^{tf(Y_z+Z_z)}$ is non-decreasing in $(Z_k)_z$, while the other function is non-increasing in $(Z_k)_z$. Next, if two of the three reals $(Z_k)_z$, $(Y_k)_z$ and $(Y'_k)_z$ have a different sign, then the function $M_p$ is constant in $(Z_k)_z$. Hence, we can use (4.7) and if we denote by $I$ the set of indices for which $(Z_k)_z$ is positive, we get

$$\mathcal{D} = \sum_{I \subset \{1,\ldots,d\}} \mathbb{E}_{Z_z}\left[Q_p\left(\frac{M_p(|u|, Y_z, Y'_z, Z_z)}{|u|d^{1/(2p-2)}}\right)e^{tf(Y_z+Z_z)}\mathbf{1}_{(Z_k)_z \geq 0 \text{ iff } \{k \in I\}}\right]$$

$$\leq \sum_{I \subset \{1,\ldots,d\}} \mathbb{E}_{Z_z}\left[Q_p\left(\frac{M_p(|u|, Y_z, Y'_z, Z_z)}{|u|d^{1/(2p-2)}}\right)\right]\mathbb{E}_{Z_z}[e^{tf(Y_z+Z_z)}\mathbf{1}_{\{(Z_k)_z \geq 0 \text{ iff } k \in I\}}]$$

$$= \mathbb{E}_{Z_z}\left[Q_p\left(\frac{M_p(|u|, Y_z, Y'_z, Z_z)}{|u|d^{1/(2p-2)}}\right)\right]\mathbb{E}_{Z_z}[e^{tf(Y_z+Z_z)}].$$

Integrating in $Y'_z$ and then in $Y_z$ leads to

$$\mathcal{D} \leq \sup_{y \in \mathbb{R}^d} \mathbb{E}_{Y'_z,Z_z}\left[Q_p\left(\frac{M_p(|u|, y, Y'_z, Z_z)}{|u|d^{1/(2p-2)}}\right)\right]\mathbb{E}[e^{tf(V)}],$$

where there is no index $z$ in the last expectation since all the marginals are the same. In this way, we have proven the second step, with

$$F(p, d, u) = \sup_{y \in \mathbb{R}^d} \mathbb{E}_{Y'_z,Z_z}\left[Q_p\left(\frac{M_p(|u|, y, Y'_z, Z_z)}{|u|d^{1/(2p-2)}}\right)\right].$$

**Proof of Step 3.** Since $Q_p$ is continuous and decreasing, let us denote by $V_p$ its reciprocal. In order to bound $F(p, d, u)$, we start by evaluating the probability

$$\mathbb{P}_{Y'_z,Z_z}\left[Q_p\left(\frac{M_p(|u|, y, Y'_z, Z_z)}{|u|d^{1/(2p-2)}}\right) \geq s\right]$$

$$= \mathbb{P}_{Y'_z,Z_z}\left[\frac{M_p(|u|, y, Y'_z, Z_z)}{|u|d^{1/(2p-2)}} \leq V_p(s)\right],$$

which is, in fact, zero if $s \geq a = Q_p(1/d^{1/(2p-2)})$.

Suppose that $y$ has $k$ positive coordinates and $d - k$ negative coordinates. Let $I_+$ be the set of $i$ such that $y_i > 0$ and let $I_-$ be the set of $i$ such that $y_i < 0$. Set

$$m_q^+(z) = \mathbb{E}(|Y_z^+ + Z_z^+|^q),$$

$$\overline{m_q}(z) = \sup(m_q^+(z), m_q^-(z)),$$

and define $m_q^-(z)$ and $\underline{m_q}(z)$ likewise, using $Y_z^-$ and $Z_z^-$. Then,

$$\mathbb{P}_{Y_z', Z_z}\left[Q_p\left(\frac{M_p(|u|, y, Y_z', Z_z)}{|u|d^{1/(2p-2)}}\right) \geq s\right]$$

$$= \mathbb{P}_{Y_z', Z_z}\left[\sum_{i \in I_+} |Y_z'^+ + Z_z^+|^p + \sum_{i \in I_-} |Y_z'^- + Z_z^-|^p \leq d^{p/(2p-2)}|u|^p V_p(s)^p - |u|^p\right]$$

$$\leq \mathbb{P}_{Y_z', Z_z}\left[\sum_{i \in I_+} |Y_z'^+ + Z_z^+|^p + \sum_{i \in I_-} |Y_z'^- + Z_z^-|^p \leq d^{p/(2p-2)}|u|^p V_p(s)^p\right].$$

Now, if $s$ is such that $d^{p/(2p-2)}|u|^p V_p(s)^p \leq d\underline{m_p}(z)/2$, which is equivalent to saying that

$$s \geq b = Q_p\left(\left[\frac{d^{(p-2)/(2p-2)} \underline{m_p}(z)}{2|u|^p}\right]^{1/p}\right), \tag{4.9}$$

then

$$\mathbb{P}_{Y_z', Z_z}\left[Q_p\left(\frac{M_p(|u|, y, Y_z', Z_z)}{|u|d^{1/(2p-2)}}\right) \geq s\right]$$

$$\leq \mathbb{P}_{Y_z', Z_z}\left[\sum_{i \in I_+} |Y_z'^+ + Z_z^+|^p + \sum_{i \in I_-} |Y_z'^- + Z_z^-|^p \leq \frac{d\underline{m_p}(z)}{2}\right]$$

$$\leq \mathbb{P}_{Y_z', Z_z}\left[\sum_{i \in I_+} |Y_z'^+ + Z_z^+|^p + \sum_{i \in I_-} |Y_z'^- + Z_z^-|^p\right.$$

$$\left. \leq km_p^+(z) + (d-k)m_p^-(z) - \frac{d\underline{m_p}(z)}{2}\right].$$

Using the Bienaymé–Chebyshev inequality, we obtain

$$\mathbb{P}_{Y_z', Z_z}\left[Q_p\left(\frac{M_p(|u|, y, Y_z', Z_z)}{|u|d^{1/(2p-2)}}\right) \geq s\right] \leq 4\frac{km_{2p}^+(z) + (d-k)m_{2p}^-(z)}{d^2 \underline{m_p}^2(z)}$$

$$\leq \frac{4\overline{m_{2p}}(z)}{d\underline{m_p}^2(z)}.$$

Integrating this last probabilistic inequality gives

$$\mathbb{E}_{Y_z', Z_z}\left[Q_p\left(\frac{M_p(|u|, y, Y_z', Z_z)}{|u|d^{1/(2p-2)}}\right)\right] \leq b + 4\left(\frac{a}{d}\right)\frac{\overline{m_{2p}}(z)}{\underline{m_p}^2(z)}.$$

Recalling that $a = Q_p(1/d^{1/(2p-2)})$, so that $a/d \leq 2^{2p-2}$, this entails

$$F(p,d,u) \leq b + 2^{2p} \frac{\overline{m_{2p}}(z)}{\overline{m_p}^2(z)}. \tag{4.10}$$

Then, using Proposition 1, where all the $b_k$ are here equal to 1, together with Step 2 and (4.10), we get Step 3.

**The positive case.** When $X$ has almost surely non-negative coordinates, the following improvements are possible. First, (4.5) can be replaced by

$$|f(X + ue_k) - f(X)| \leq p|u| \left( \frac{|X_k| + |u|}{(2\|X\|_p^p + |u|^p)^{1/p}} \right)^{p-1}. \tag{4.11}$$

Indeed, we can use the proof of (4.5), but here we only have to consider the case $c_k \geq 0$.

Second, if $X$ has non-negative coordinates, we can apply (4.7) directly, without introducing $M_p$. In this way, (4.10) can be replaced by

$$F(d,p,u) \leq b + 4 \left( 1 + \frac{1}{2^{1/p}} \right)^{2p-2} \frac{m_{2p}}{m_p^2}, \tag{4.12}$$

where $m_p$ and $m_{2p}$ now stand for the usual moments, that is, $m_q = \mathbb{E}[X^q]$ for $q = p, 2p$. $\square$

## 5. Proof of Theorem 3

Let

$$g(x) = \frac{x}{R} \log^+ \left( \frac{\lambda x}{R} \right).$$

Note that $g$ is a bijection from $[R/\lambda, +\infty[$ to $[0, +\infty[$. To prove Theorem 3, it is thus sufficient to determine for which $\lambda$, $\mathbb{E}e^{g(\|X\|_2)}$ converges. Let $\varepsilon > 0$ and let $c = \max((1+\varepsilon)\mathbb{E}\|X\|_2, R/\lambda)$. Then,

$$\mathbb{E}e^{g(\|X\|_2)} = \int_0^{+\infty} \mathbb{P}(e^{g(\|X\|_2)} \geq t) \, dt \leq e^{g(c)} + \int_{e^{g(c)}}^{+\infty} \mathbb{P}(e^{g(\|X\|_2)} \geq t) \, dt.$$

Setting $t = e^{g(c+u)}$ in the last integral and applying Corollary 2, we get

$$\mathbb{E}e^{g(\|X\|_2)} \leq e^{g(c)} + \int_0^{+\infty} e^{u/R - (u/R + V_\varepsilon^2/R^2)\log(1 + Ru/V_\varepsilon^2)}$$
$$\times \left[ \frac{1}{R} + \frac{1}{R} \log^+ \left( \frac{\lambda(c+u)}{R} \right) \right] e^{(c+u)/R \log^+ (\lambda(c+u)/R)} \, du.$$

For this last integral to converge, it is sufficient that the power of $e^{u/R}$ be negative, that is, that

$$1 + \log^+\left(\frac{\lambda(c+u)}{R}\right) - \log\left(1 + \frac{Ru}{V_\varepsilon^2}\right) < 0,$$

at least for $u$ large enough. But, if $u$ tends to $+\infty$ in the expression above, we obtain that

$$1 + \log\left(\frac{\lambda V_\varepsilon^2}{R^2}\right) < 0, \quad \text{that is,} \quad \lambda < \frac{R^2 e^{-1}}{V_\varepsilon^2}.$$

Next, we need to pass from $V_\varepsilon$ to $V$. But, if $\lambda < R^2 e^{-1}/V^2$, then there exists $\varepsilon$ (quite large, of course) such that $\lambda < R^2 e^{-1}/V_\varepsilon^2$ and this concludes the proof.

## Acknowledgements

Many thanks to the referees for their detailed comments which greatly helped the readability and presentation of this paper. This research was done, in part, while the first author visited L'École Normale Supérieure whose hospitality and support are gratefully acknowledged. The work of the second author was supported by the French Agence Nationale de la Recherche (ANR), grant ATLAS (JCJC06 137446) "From Applications to Theory in Learning and Adaptive Statistics". The work of the third author was supported by the French Agence Nationale de la Recherche (ANR), grant ATLAS (JCJC06 137446) "From Applications to Theory in Learning and Adaptive Statistics".

## References

[1] Bobkov, S.G. and Houdré, C. (2000). Weak dimension-free concentration of measure. *Bernoulli* **6** 621–632. MR1777687
[2] Bobkov, S.G. and Ledoux, M. (1997). Poincaré inequalities and Talagrand's concentration phenomenon for the exponential measure. *Probab. Theory Related Fields* **107** 383–400. MR1440138
[3] Borovkov, A.A. and Utev, S.A. (1983). An inequality and a characterization of the normal distribution connected with it. *Theory Probab. Appl.* **28** 209–218. MR0700206
[4] Donoho, D., Johnstone, I., Kerkyacharian, G. and Picard, D. (1996). Density estimation by wavelet thresholding. *Ann. Statist.* **24** 508–539. MR1394974
[5] Houdré, C., Pérez-Abreu, V. and Surgailis, D. (1998). Interpolation, correlation identities and inequalities for infinitely divisible variables. *J. Fourier Anal. Appl.* **4** 651–668. MR1665993
[6] Houdré, C. (2002). Remarks on deviation inequalities for functions of infinitely divisible random vectors. *Ann. Probab.* **30** 1223–1237. MR1920106
[7] Houdré, C. and Reynaud-Bouret, P. (2004). Concentration for infinitely divisible vectors with independent components. Available at http://arXiv.org/abs/math.PR/0606752.
[8] Massart, P. (2007). *Concentration Inequalities and Model Selection. Ecole d'Été de Probabilités de Saint-Flour 2003. Lecture Notes in Math.* **1896**. Berlin: Springer. MR2319879

 [9] Maurey, B. (1991). Some deviation inequalities. *Geom. Funct. Anal.* **1** 188–197.
[10] Reynaud-Bouret, P. (2003). Adaptive estimation of the intensity of inhomogeneous Poisson processes via concentration inequalities. *Probab. Theory Related Fields* **126** 103–153. MR1981635
[11] Rosiński, J. (1996). Remarks on strong exponential integrability of vector-valued random series and triangular arrays. *Ann. Probab.* **23** 464–473. MR1330779
[12] Sato, K.-I. (1999). *Lévy Processes and Infinitely Divisible Distributions*. Translated from the 1990 Japanese original. Revised by the author. *Cambridge Studies in Advanced Mathematics* **68**. Cambridge Univ. Press. MR1739520
[13] Talagrand, M. (1991). A new isoperimetric inequality for product measure, and the concentration of measure phenomenon. *Israel Seminar* (*GAFA*). *Lecture Notes in Math.* **1469** 91–124. Berlin: Springer. MR1122615